\documentclass[a4paper,12pt]{amsart}
\usepackage{amssymb}
\usepackage{ifthen}
\usepackage{graphicx}
\usepackage{float}
\usepackage{caption}
\usepackage{subcaption}
\usepackage{cite}
\usepackage{amsfonts}
\usepackage{amscd}
\usepackage{amsxtra}
\usepackage{color}

\newtheorem{thm}{Theorem}
\newtheorem{cor}{Corollary}
\newtheorem{lem}{Lemma}
\newtheorem{rem}{Remark}

\newtheorem{conj}{Conjecture}
\newtheorem{prob}{Problem}

\theoremstyle{definition}
\newtheorem{defn}{Definition}[section]
\newtheorem{example}{Example}

\newenvironment{pf}[1][]{%
 \vskip 1mm
 \noindent
 \ifthenelse{\equal{#1}{}}%
  {{\slshape Proof. }}%
  {{\slshape #1.} }%
 }%
{\qed\bigskip}

\newcounter{alphabet}

\newenvironment{Thm}[1][]{\refstepcounter{alphabet}%
\bigskip%
\noindent%
{\bf Theorem \Alph{alphabet}}%
\ifthenelse{\equal{#1}{}}{}{ (#1)}%
{\bf .} \itshape}{\vskip 8pt}
\newenvironment{Lem}[1][]{\refstepcounter{alphabet}%
\bigskip%
\noindent%
{\bf Lemma \Alph{alphabet}}%
{\bf .} \itshape}{\vskip 8pt}

\newcommand{\B}{{\mathcal B}}

\newcommand{\ID}{{\mathbb D}}

\newcommand{\dist}{{\operatorname{dist}}}

\def\be{\begin{equation}}
\def\ee{\end{equation}}

\newcommand{\bee}{\begin{enumerate}}
\newcommand{\eee}{\end{enumerate}}

\newcommand{\blem}{\begin{lem}}
\newcommand{\elem}{\end{lem}}
\newcommand{\bthm}{\begin{thm}}
\newcommand{\ethm}{\end{thm}}
\newcommand{\bcor}{\begin{cor}}
\newcommand{\ecor}{\end{cor}}
\newcommand{\beg}{\begin{example}}
\newcommand{\eeg}{\end{example}}
\newcommand{\begs}{\begin{examples}}
\newcommand{\eegs}{\end{examples}}
\newcommand{\bdefe}{\begin{defn}}
\newcommand{\edefe}{\end{defn}}
\newcommand{\bprob}{\begin{prob}}
\newcommand{\eprob}{\end{prob}}
\newcommand{\bques}{\begin{ques}}
\newcommand{\eques}{\end{ques}}
\newcommand{\bei}{\begin{itemize}}
\newcommand{\eei}{\end{itemize}}
\newcommand{\bcon}{\begin{conj}}
\newcommand{\econ}{\end{conj}}
\newcommand{\bcons}{\begin{conjs}}
\newcommand{\econs}{\end{conjs}}
\newcommand{\bprop}{\begin{propo}}
\newcommand{\eprop}{\end{propo}}
\newcommand{\br}{\begin{rem}}
\newcommand{\er}{\end{rem}}
\newcommand{\brs}{\begin{rems}}
\newcommand{\ers}{\end{rems}}
\newcommand{\bo}{\begin{obser}}
\newcommand{\eo}{\end{obser}}
\newcommand{\bos}{\begin{obsers}}
\newcommand{\eos}{\end{obsers}}
\newcommand{\bpf}{\begin{pf}}
\newcommand{\epf}{\end{pf}}
\newcommand{\ba}{\begin{array}}
\newcommand{\ea}{\end{array}}
\newcommand{\beq}{\begin{eqnarray}}
\newcommand{\beqq}{\begin{eqnarray*}}
\newcommand{\eeq}{\end{eqnarray}}
\newcommand{\eeqq}{\end{eqnarray*}}

\newcommand{\ds}{\displaystyle}

\newcounter{minutes}
\divide\time by 60
\newcounter{hours}
\multiply\time by 60 \addtocounter{minutes}{-\time}
\begin{document}

\bibliographystyle{amsplain}


\title[Note on Weighted Bohr's Inequality]{Note on Weighted Bohr's Inequality}

\thanks{
File:~\jobname .tex,
          printed: \number\day-\number\month-\number\year,
          \thehours.\ifnum\theminutes<10{0}\fi\theminutes}


\author[R. Vijayakumar]{Ramakrishnan Vijayakumar}
\address{
R. Vijayakumar, Department of Mathematics,
Indian Institute of Technology Madras, Chennai-600 036, India.
}
\email{mathesvijay8@gmail.com}


\begin{abstract}
In this paper, first we give a new generalization of the Bohr's inequality for the class of bounded analytic functions $\mathcal{B'}$ and for the class of sense-preserving $K$-quasiconformal harmonic mappings of the form $f=h+\overline{g},$ where $h \in \mathcal{B'}.$ Finally we give a new generalization of the Bohr's inequality for the class of analytic functions subordinate to univalent functions and for the class of sense-preserving $K$-quasiconformal harmonic mappings of the form $f=h+\overline{g},$ where $h$ is subordinated to some analytic function.
\end{abstract}

\maketitle
\pagestyle{myheadings}
\markboth{R. Vijayakumar}{Weighted Bohr's Inequality}

\section{Introduction and Preliminaries}
Throughout we let $\mathcal B$ denote the class of all analytic functions $\omega$ in the open unit disk $\mathbb{D}=\{z\in \mathbb{C}:\,|z|<1\}$
such that $|\omega(z)|\leq 1$ for all $z \in \ID.$ Bohr's inequality says that if $f \in \mathcal B$ and $f(z)=\sum_{n=0}^{\infty} a_n z^n,$ then we have
$$ \sum_{n=0}^{\infty} |a_{n}| r^n \leq 1
$$
for all $z \in \ID$ with $|z|=r \leq \frac{1}{3}.$ This inequality was discovered by Bohr in 1914 \cite{Bohr}. Bohr actually obtained the inequality
for $|z| \leq \frac{1}{6}.$ Later M. Riesz, I. Schur and F. W. Wiener independently, established the inequality for $|z| \leq \frac{1}{3}$ and showed that
$\frac{1}{3}$ is sharp. The number is  $\frac{1}{3}$  is called Bohr radius
for the family $\mathcal B$. A space of analytic or harmonic functions $f$ in $\ID$ is said to have
Bohr's phenomenon if an inequality of this type holds in some disk of radius $\rho>0$ and for all such functions in unit ball of the space. In \cite{BenDahKha},
it is shown that not every space of functions has Bohr's phenomenon. On the other hand,
Abu-Muhanna \cite{Abu-1} proved the existence of Bohr phenomenon in the case of subordination and bounded harmonic classes.
Many mathematicians have contributed towards the understanding of this problem in several settings \cite{Bomb-1962,BombBour-2004}.
Extensions of Bohr's inequality to more general domains or higher dimensional spaces were investigated by many. See for instance, \cite{BoasKhavin-97-4,DjaRaman-2000,HHK2009}.
We refer to the recent survey on this topic by Abu-Muhanna et al. \cite{AAPon1} and Garcia et al. \cite{GarMasRoss-2018}, for the importance and the several other results.
For certain recent results, see \cite{AlkKayPon1,KayPon1,KayPon2,KayPon3}.

More generally, a harmonic version of Bohr's inequality was discussed by Kayumov et al. \cite{KPS}. For certain other results on harmonic  Bohr's inequality, we refer to \cite{EvPoRa, KPS}. Recently, a new generalization of Bohr's ideas was introduced and investigated by Kayumov et al. \cite{KaykhaPo}. In order to make the statement of the recent generalization,
we need to introduce some basic notations.

Let ${\mathcal F}$ denote the set of all sequences  $\{\varphi_{n}(r)\}_{n=0}^{\infty}$ of nonnegative continuous functions in $[0, 1)$ such that the series $\sum_{n=0}^{\infty} \varphi_{n}(r)$ converges locally uniformly with respect to $r \in [0, 1)$.  Let ${\mathcal F}_{dec} \subset {\mathcal F}$ consist of decreasing sequences
of functions from ${\mathcal F}$, and for convenience, we let $\Phi_1(r)=\sum_{n=1}^{\infty} \varphi_{n}(r)$ so that $\Phi_1'(r)=\sum_{n=1}^{\infty} \varphi_{n}'(r)$
whenever each $\varphi_{n}$ $(n\ge 1)$ is differentiable on $[0,1]$.

\begin{Thm} \label{KKP16-th1}
{\rm (\cite{KaykhaPo})}
Let $f\in \mathcal{B}$, $f(z)=\sum_{k=0}^\infty a_k z^k$, and $p\in (0,2]$. If
$ \varphi_0(r) > (2/p) \Phi_1(r),$ then  the following sharp inequality holds:
$$
B_{f}(\varphi, p, r):=|a_0|^p \varphi_0(r) + \sum_{k=1}^{\infty} |a_k| \varphi_k(r) \leq \varphi_0(r) ~\mbox{ for all $r \leq R$},
$$
where $R$ is the minimal positive root of the equation
$\varphi_0(x) = (2/p) \Phi_1(x). $
In the case when $\varphi_0(x) < (2/p)\Phi_1(x) $
in some interval $(R,R+\varepsilon)$, the number $R$ cannot be improved. If the functions $\varphi_k(x)$ ($k\geq 0$) are smooth functions, then the last condition
is equivalent to the inequality $\varphi_0'(R) <(2/p)\Phi_1'(R). $
\end{Thm}

Further investigation and refinements of several earlier known results on Bohr-type inequality, we refer to \cite{PoViWirth-3}.


For two analytic functions $f$ and $g$ in $\ID$, we say that $g$ is subordinate to $f$ (denoted simply by $g \prec f$)  if there exists a function $\omega,$ analytic in $\ID$
with $\omega(0)=0$ and $|\omega(z)|<1,$ satisfying $g=f \circ \omega.$
We denote the class of all analytic functions $g$ in $\ID$ that are subordinate to a fixed function $f$ by $\mathcal{S}(f)$, and $f(\mathbb{D})=\Omega.$
We say that $\mathcal{S}(f)$ has Bohr's phenomenon if for any $g(z)= \sum_{n=0}^{\infty}b_{n}z^{n} \in \mathcal{S}(f)$ and
$f(z)=\sum_{n=0}^{\infty}a_{n}z^{n}$, there is a $\rho_{0}$, $0 < \rho_{0} \leq 1$, so that
$$\sum_{n=1}^{\infty}|b_{n}z^{n}|\leq \dist(f(0), \partial \Omega),
$$
for $|z|< \rho_{0}.$ We remark that the class $\mathcal{S}(f)$ has Bohr's phenomenon when $f$ is univalent (see \cite[Theorem 1]{Abu-1}). For each
$f(z)=\sum_{k=0}^\infty a_k z^k$ belonging to $\mathcal{B}$, it is well-known that $ |a_{n}|\leq 1-|a_{0}|^2$ for all $n\geq 1$.
Besides the fact that  $1-|a_0|\le 1-|a_0|^2$ for $|a_0|\le 1$, as demonstrated by Aizenberg and Vidras (see \cite[p. 736]{AizenbergVidras},
there exists a nice subclass of functions $f\in\mathcal{B}$ for which $ |a_{n}|\leq 1-|a_{0}|$ all $n\geq 1$. We now recall this result.

\begin{Thm}\label{AizVed}{\rm (\cite{AizenbergVidras})}
Let $f\in \mathcal B,$ such that the Taylor coefficients $a_{mn} =0$ for a given $m>1$ and all $n \geq 1.$ Then
$|a_n|\le 1-|a_0|$ for all  $n\ge1.$
\end{Thm}

Thus, it is natural to consider
$$\mathcal{B'}=\left \{f(z)=\sum_{k=0}^\infty a_k z^k\in \mathcal{B}:\,  |a_{n}|\leq 1-|a_{0}| ~\mbox{ for all $n \geq 1$} \right \}.
$$
In \cite[Theorem 1]{Alkha}, it was shown that  the Bohr radius for functions in $\B'$ is $\frac12$, and the constant $1/2$ cannot be improved.

In this article, we first investigate the Bohr radius for the family $\B'$ in a general setting which is indeed an analog of Theorem~A 
for the family $\mathcal{B'}$ (See Theorem \ref{thm 1}). Our second result (Theorem \ref{thm 2}) extends  Theorem~A 
to the case of sense-preserving $K$-quasiconformal harmonic mappings of the form $f=h+\overline{g},$ where $h \in \mathcal{B'}.$
In  Section \ref{PV8-sec4}, we establish that the family $\mathcal{S}(f)$ has Bohr's phenomenon in our new setting (see Theorems \ref{thm 3} and \ref{thm 4}), especially
when $f$ is either univalent or convex (univalent) in $\ID$. Finally, we extend this result  (Theorem \ref{thm 5}) for sense-preserving $K$-quasiconformal harmonic mappings.


\section{Bohr radius for a special family of analytic functions}

The following theorem displays the sharp Bohr radius for $\mathcal{B'}$.

\begin{thm}\label{thm 1}
Let $f\in \mathcal{B'}$,   $f(z)=\sum_{n=0}^{\infty} a_n z^n$,  and $p\in (0,1].$ If $\{\varphi_{n}(r)\}_{n=0}^{\infty} \in {\mathcal F}$
such that $\Phi_1(r)=\sum_{n=1}^{\infty} \varphi_{n}(r)$, and satisfies the inequality
\be\label{eq 1}
\varphi_{0}(r) \geq  \frac{1}{p}\Phi_1(r).
\ee
 Then the following sharp inequality holds:
\be\label{eq 2}
B_{f}(\varphi, p, r):= |a_{0}|^p \varphi_{0}(r)+ \sum_{n=1}^{\infty} |a_{n}| \varphi_{n}(r) \leq \varphi_{0}(r)~\mbox{for all}~ r \leq R,
\ee
where $R$ is the minimal positive root of the equation
\beqq
\varphi_{0}(x) =  \frac{1}{p}\Phi_1(x).
\eeqq
In the case when $\varphi_{0}(x) <  \frac{1}{p}\Phi_1(x)$ in some interval $(R, R+\epsilon),$ the number $R$ cannot be improved.
	
\end{thm}
\begin{pf}
Let $f\in \mathcal{B'}.$ Then $|a_{n}|\leq 1-|a_{0}| ~\mbox{for all} \ \ n \geq 1$ and thus, 
we get that
\beq
|a_{0}|^p \varphi_{0}(r)+ \sum_{n=1}^{\infty} |a_{n}| \varphi_{n}(r)
&\leq& |a_{0}|^p \varphi_{0}(r)+ (1-|a_{0}|) \Phi_1(r) \nonumber \\
&=&\varphi_{0}(r) +(1-|a_{0}|) \left[\Phi_1(r)-\left(\frac{1-|a_{0}|^p}{1-|a_{0}|}\right)\varphi_{0}(r) \right] \nonumber \\
&\leq& \varphi_{0}(r) +(1-|a_{0}|) \left[\Phi_1(r) -p \varphi_{0}(r) \right] \nonumber \\
&\leq& \varphi_{0}(r), \ \ \text{by Eqn. \eqref{eq 1}} \nonumber,
\eeq	
for all $r \leq R,$ by the definition of $R.$ In the third inequality above, we have used the fact that the function
$$B(x)=\frac{1-x^p}{1-x}, \quad x\in [0,1),
$$
is decreasing on $[0,1)$ for $0<p \leq 1$ so that
$$ B(x) \geq \lim_{x\rightarrow 1^{-}}\frac{1-x^p}{1-x}\, =p.
$$

This proves the desired inequality \eqref{eq 2}. Now let us prove that $R$ is an optimal number. For $a\in [0,1),$ we consider the function	
$$ f(z) =\frac{a-(1-a+a^2)z}{1-a z}=a-(1-a)\sum\limits_{n=1}^\infty a^{n-1}z^n,  ~ z \in \mathbb{D}.
$$	
A simple exercise shows that $f \in \mathcal{B'}.$ For this function, we have
\beqq
|a_{0}|^p \varphi_{0}(r)+ \sum_{n=1}^{\infty} |a_{n}| \varphi_{n}(r)&=& a^p\varphi_{0}(r) +(1-a) \sum_{n=1}^{\infty} a^{n-1} \varphi_{n}(r) \\
& =& \varphi_{0}(r)+p(1-a)\left[\frac{1}{p} \sum_{n=1}^{\infty}a^{n-1}\varphi_{n}(r)- \varphi_{0}(r)\right]\\
&& \quad +(1-a)\left[\left(p- \frac{1-a^p}{1-a}\right)\varphi_{0}(r)\right].\\
\eeqq
 Now it is easy to see that number is $> \varphi_{0}(r)$ when $a$ is close to 1. The proof of the theorem is complete.
\end{pf}	

\br Note that the function $B(x)$ in the above proof is increasing on $[0,1)$ for $p \geq 1$ so that $B(x) \geq B(0)=1.$ This means that  the inequality
\eqref{eq 2} holds for $r \leq \frac{1}{2}$ in the case when $\varphi_{n}(r)=r^{n} \, (n \geq 1).$
\er

\begin{cor}\label{cor 1}
Suppose that $f \in \mathcal{B'},$ $f(z)= \sum_{n=0}^{\infty} a_{n}z^{n}$, and $p\in (0,1].$ Then
\beqq
|a_0|^p+ \sum_{n=1}^\infty |a_n|r^n \leq 1 ~\mbox{for}~ r \leq R(p)=\frac{p}{1+p},
\eeqq
and the constant $R(p)$ cannot be improved.
\end{cor}

 The case $p=1$ of Corollary \ref{cor 1} is the Bohr inequality for special family of bounded analytic functions $\mathcal{B'},$ obtained in \cite[Theorem 1]{Alkha}.
%

\section{Bohr radius for harmonic mappings as an extension of Theorem \ref{thm 1}}

We recall that a sense-preserving harmonic mappings $f$ of the form $f=h+\overline{g},$
is said to be $K$-quasiconformal if $|g'(z)|\leq k|h'(z)|$ in the unit disk, for $k=\frac{K-1}{K+1} \in [0,1].$
See \cite{KPS} for discussion on Bohr radius for quasiconformal mappings.

\begin{Lem}\label{lem 1}
{\rm (\cite{PoViWirth-3})}
Let $\{\psi_{n}(r)\}_{n=1}^{\infty}$ be a decreasing sequence of nonnegative   functions in $[0,r_\psi)$, and
$g, h$ be analytic functions in the unit disk $\ID$ such that $|g'(z)|\leq k |h'(z)|$ in $\ID$ and for some $k \in [0,1],$ where
$h(z)=\sum_{n=0}^\infty a_nz^n$ and $g(z)=\sum_{n=0}^\infty b_nz^n$. Then
$$\sum_{n=1}^\infty |b_n|^2 \psi_{n}(r) \leq k^2 \sum_{n=1}^\infty |a_n|^2 \psi_{n}(r) \ \ \text{for}\ \ r\in [0,r_\psi).
$$
\end{Lem}

Next, we find Bohr radius for the family of sense-preserving $K$-quasiconformal harmonic mappings of the form $f=h+\overline{g},$ where $h \in \mathcal{B'}$
and show the sharpness of it.

\begin{thm}\label{thm 2}
Suppose that $f(z)= h(z)+\overline{g(z)}=\sum_{n=0}^\infty a_nz^n +\overline{\sum_{n=1}^\infty b_nz^n}$ is harmonic mapping of the disk $\ID$
such that $|g'(z)|\leq k |h'(z)|$ in $\ID$ and for some $k \in [0,1],$ where $h \in \mathcal{B'}.$ Assume that $\varphi_{0}(r)=1$ and
$\{\varphi_{n}(r)\}_{n=0}^{\infty}$ belongs to ${\mathcal F}_{dec}$ with $\Phi_1(r)=\sum_{n=1}^{\infty} \varphi_{n}(r)$,
and $p\in (0,1].$ If
\beq\label{eq 3}
p \geq (1+k)\Phi_1(r),
\eeq
then the following sharp inequality holds:
\beq\label{eq 4}
|a_{0}|^p+\sum_{n=1}^\infty |a_n| \varphi_{n}(r)+\sum_{n=1}^\infty |b_n| \varphi_{n}(r) \leq \|h\|_{\infty} \  \ \text{for all}\  \ r \le R,
\eeq
where $R$ is the minimal positive root of the equation
$$ p=(1+k)\Phi_1(x).
$$
In the case when $p<(1+k)\Phi_1(x)$	in some interval $(R, R+\epsilon),$ the number $R$ cannot be improved.
\end{thm}	
\begin{pf}
For simplicity, we suppose that $\|h\|_{\infty}=1.$ For $h\in \mathcal{B'},$ gives the inequality $|a_{n}|\leq 1-|a_{0}| ~\mbox{for all} \ \ n \geq 1.$ By assumption $|g'(z)|\leq k |h'(z)|$ in $\ID,$ where $k \in[0,1]$
and so, by Lemma~C 
it follows that
\beqq	
\sum_{n=1}^\infty |b_n|^2 \varphi_{n}(r) \leq k^2 \sum_{n=1}^\infty |a_n|^2 \varphi_{n}(r) \leq k^2 (1-|a_{0}|)^2 \sum_{n=1}^{\infty} \varphi_{n}(r) =k^2 (1-|a_{0}|)^2\Phi_1(r).
\eeqq	
Consequently, it follows from the classical Schwarz inequality that
\beqq	
\sum_{n=1}^\infty |b_n| \varphi_{n}(r) \leq  \sqrt{\sum_{n=1}^\infty |b_n|^2 \varphi_{n}(r)} \sqrt{\sum_{n=1}^{\infty}\varphi_{n}(r)}
\leq  
k (1-|a_{0}|)\Phi_1(r)
\eeqq
and thus, as in the proof of Theorem \ref{thm 1}, we get that
\beq
|a_{0}|^p+\sum_{n=1}^\infty |a_n| \varphi_{n}(r)+\sum_{n=1}^\infty |b_n| \varphi_{n}(r) &\leq& |a_{0}|^p + (1-|a_{0}|)(1+k) \Phi_1(r) \nonumber \\
&=&1 +(1-|a_{0}|) \left[(1+k)\Phi_1(r)-\left(\frac{1-|a_{0}|^p}{1-|a_{0}|}\right) \right] \nonumber \\
&\leq & 1 +(1-|a_{0}|) \left[(1+k)\Phi_1(r)-p  \right] \nonumber \\
&\leq& 1, \ \ \text{by Eqn. \eqref{eq 3},} \nonumber
		\eeq	
		for all $r \leq R,$ by the definition of $R.$ This proves the desired inequality \eqref{eq 4}.
Now let us prove that $R$ is an optimal number. We consider the function
$$ h(z) =\frac{a-(1-a+a^2)z}{1-a z}=a-(1-a)\sum\limits_{n=1}^\infty a^{n-1}z^n, a\in [0,1), z \in \mathbb{D}
$$	
and $g(z)= \lambda k h(z),$ where $|\lambda|=1$. Then it is a simple exercise to see that

\vspace{8pt}
$\ds |a_{0}|^p+\sum_{n=1}^\infty |a_n| \varphi_{n}(r)+\sum_{n=1}^\infty |b_n| \varphi_{n}(r)
$
\beqq
 &=& a^{p}+(1-a)\sum\limits_{n=1}^\infty a^{n-1}\varphi_{n}(r)+k(1-a)\sum\limits_{n=1}^\infty a^{n-1}\varphi_{n}(r)\\
 & =& 1+p(1-a)\left[\frac{1}{p}(1+k) \sum_{n=1}^{\infty}a^{n-1}\varphi_{n}(r)-1 \right] +(1-a)\left(p- \frac{1-a^p}{1-a}\right).
\eeqq
Now it is easy to see that number is $> 1$ when $a$ is close to 1. The proof of the theorem is complete.
\end{pf}
	
\begin{cor}
Suppose that $f(z)= h(z)+\overline{g(z)}=\sum_{n=0}^\infty a_nz^n +\overline{\sum_{n=1}^\infty b_nz^n}$ is a sense-preserving $K$-quasiconformal
harmonic mapping of the disk $\ID$, i.e. $|g'(z)|\leq k |h'(z)|$ in $\ID$ for some $k=\frac{K-1}{K+1} \in[0,1],$ where $h \in \mathcal{B'}$.
Then we have the sharp inequality
\beq \label{eq 5}
|a_{0}|^p+\sum_{n=1}^\infty |a_n|r^n+\sum_{n=1}^\infty |b_n|r^n \leq 1 \ \text{for} \ r \le R_{k}(p)
\eeq
where $p\in (0,1]$, and
$$R_{k}(p)=\frac{p}{k+1+p}=\frac{p(K+1)}{(p+2)K+p}
$$
and the constant $R_{k}(p)$ cannot be improved.
\end{cor}	

In particular, the case $p=1$ in \eqref{eq 5} yields the recently  obtained result \cite[Theorem 2]{Alkha}.

\section{Bohr phenomenon in subordination}\label{PV8-sec4}

The following lemma will be used to prove that the family $\mathcal{S}(f)$ has Bohr's phenomenon in our new setting (see Theorem \ref{thm 3}).

\begin{Lem}\label{lem A}
{\rm \cite[p. 195-196]{Duren}}
Let $f$ be an analytic univalent map from $\mathbb{D}$ onto a simply connected domain $\Omega:=f(\mathbb{D})$ and $g(z)=\sum_{n=0}^{\infty}b_{n}z^{n} \prec f(z).$ Then
$$ \frac{1}{4}|f'(0)| \leq \dist(f(0), \partial \Omega) \leq |f'(0)|,\ \mbox{and}~|b_{n}|\leq n|f'(0)|\leq 4n\ \dist(f(0), \partial \Omega).
$$
\end{Lem}

\begin{thm}\label{thm 3}
Suppose that $g(z)=\sum_{n=0}^\infty b_nz^n \in \mathcal{S}(f)$ and $f(z)= \sum_{n=0}^\infty a_nz^n$ is univalent in $\ID.$ If $\{\varphi_{n}(r)\}_{n=1}^{\infty} \in {\mathcal F}$ satisfies the inequality
\beq\label{eqn 6}
1 \geq 4\Psi_1(r),
\eeq
where $\Psi_1(r)=\sum_{n=1}^{\infty} n\varphi_{n}(r)$, then the following sharp inequality holds:
\beq\label{eqn 7}
\sum_{n=1}^\infty |b_n| \varphi_{n}(r) \leq \dist(f(0), \partial \Omega) \  \ \text{for all}\  \ r \le R,
\eeq
where $R$ is the minimal positive root of the equation $ 1 =4\Psi_1(x). $
In the case when $1 <4\Psi_1(x)$	in some interval $(R, R+\epsilon),$ the number $R$ cannot be improved.
\end{thm}
\begin{pf}
By assumption $g\prec f$ and $f$ is univalent in $\ID.$ Then, by Lemma~D, 
we have
$$ |b_{n}| \leq 4n\ \dist(f(0), \partial \Omega).
$$
Thus, we have
\beqq
\sum_{n=1}^\infty |b_n| \varphi_{n}(r) &\leq& 4 \dist(f(0), \partial \Omega) \sum_{n=1}^\infty n \varphi_{n}(r)= 4 \dist(f(0), \partial \Omega)\Psi_1(r)\\
&\leq& \dist(f(0), \partial \Omega), \ \ \text{by Eqn.\eqref{eqn 6}},
\eeqq
for all $r \leq R,$ by the definition of $R.$ This proves the desired inequality \eqref{eqn 7}.
Now let us prove that $R$ is an optimal number. We consider the function
$$ g(z)=f(z) =\frac{z}{(1-z)^2}=\sum_{n=1}^{\infty}nz^{n}, \ \ z \in \mathbb{D}.
$$	
Then it is easy to show that
$$	\dist(f(0), \partial \Omega)=\frac{1}{4}~\mbox{and}~\sum_{n=1}^\infty |b_n| \varphi_{n}(r)=\sum_{n=1}^{\infty}n\varphi_{n}(r).
$$
Now it is easy to see that number is $> \frac{1}{4}$ when $r>R$. The proof of the theorem is complete.
\end{pf}

\br
It is a simple exercise to see that if $\varphi_{n}(r)= r^{n} \,(n\geq 1)$, then  Theorem  \ref{thm 3} yields the result of Abu-Muhanna \cite[Theorem 1]{Abu-1} with $R=3-\sqrt{8}$.
\er

The next lemma will be used to prove Theorems \ref{thm 4} and \ref{thm 5}.

\begin{Lem}\label{lem B}
{\rm \cite[p. 195-196]{Duren}}
Let $\psi$ be an analytic univalent map from $\mathbb{D}$ onto a convex domain $\Omega:=\psi(\mathbb{D})$ and $g(z)=\sum_{n=0}^{\infty}b_{n}z^{n} \prec \psi(z).$ Then
$$ \frac{1}{2}|\psi '(0)| \leq \dist(\psi(0), \partial \Omega) \leq |\psi'(0)|,\ \mbox{and}\ |b_{n}|\leq |\psi '(0)|\leq 2\ \dist(\psi(0), \partial \Omega).
$$
\end{Lem}

\begin{thm}\label{thm 4}
Suppose that $g(z)=\sum_{n=0}^\infty b_nz^n \in \mathcal{S}(f)$ and $f(z)= \sum_{n=0}^\infty a_nz^n$ is univalent and convex in $\ID.$
If $\{\varphi_{n}(r)\}_{n=0}^{\infty} \in {\mathcal F}$ satisfies the inequality
\beqq
1 \geq 2\Phi_1(r),
\eeqq
where $\Phi_1(r)=\sum_{n=1}^{\infty} \varphi_{n}(r)$, then the following sharp inequality holds:
\beqq
\sum_{n=1}^\infty |b_n| \varphi_{n}(r) \leq \dist(f(0), \partial \Omega) \  \ \text{for all}\  \ r \le R,
\eeqq
where $R$ is the minimal positive root of the equation $1 =2\Phi_1(x).$
In the case when $1 <2\Phi_1(x)$ in some interval $(R, R+\epsilon),$ the number $R$ cannot be improved.
\end{thm}
\begin{pf}
 The proof follows if we use the method of proof of Theorem \ref{thm 3} and use  Lemma~E  
 in place of by Lemma~D. 
 Sharpness follows by considering the following function
$$g(z)=f(z) =\frac{1}{1-z}=\sum_{n=0}^{\infty}z^{n} \ \mbox{for}\ z \in \mathbb{D},
$$	
so that
$$\dist(f(0), \partial \Omega)=\frac{1}{2}\ \mbox{and}\, \sum_{n=1}^{\infty} |b_n| \varphi_{n}(r)=\sum_{n=1}^{\infty}\varphi_{n}(r).
$$
Now it is easy to see that number is $> \frac{1}{2}$ when $r>R$. The proof of the theorem is complete.
\end{pf}

\br
It is a simple exercise to see that if $\varphi_{n}(r)= r^{n} \,(n\geq 1)$, then Theorem \ref{thm 4} yields the
remark of Abu-Muhanna \cite[Remark 1]{Abu-1} with $R=1/3$.
\er

\begin{thm}\label{thm 5}
Suppose that $f(z)= h(z)+\overline{g(z)}=\sum_{n=0}^\infty a_nz^n +\overline{\sum_{n=1}^\infty b_nz^n}$ is harmonic mapping of the disk $\ID$
such that $|g'(z)|\leq k |h'(z)|$ in $\ID$ and for some $k \in [0,1]$ and $h \prec \psi,$ where $\psi$ is univalent and convex in $\mathbb{D}.$  Assume that
$\{\varphi_{n}(r)\}_{n=0}^{\infty}$ belongs to ${\mathcal F}_{dec}$ and $\Phi_1(r) =\sum_{n=1}^{\infty} \varphi_{n}(r)$. If
\beq\label{eqn 8}
1 >2(1+k)\Phi_1(r),
\eeq
then the following sharp inequality holds:
\beq\label{eqn 9}
\sum_{n=1}^\infty |a_n| \varphi_{n}(r)+\sum_{n=1}^\infty |b_n| \varphi_{n}(r) \leq \dist(\psi(0), \partial \psi(\mathbb{D})) \  \ \text{for all}\  \ r \le R,
\eeq
where $R$ is the minimal positive root of the equation
$ 1=2(1+k)\Phi_1(x).
$
In the case when $1<2(1+k)\Phi_1(x)$	in some interval $(R, R+\epsilon),$ the number $R$ cannot be improved.
\end{thm}
\begin{pf}
By assumption $h \prec \psi$ and $\psi(\mathbb{D})$ is a convex domain. Then, by Lemma~E, 
we have
$$ |a_{n}| \leq 2\ \dist(\psi(0), \partial \psi(\mathbb{D})).
$$
Consequently,
$$
\sum_{n=1}^\infty |a_n| \varphi_{n}(r) \leq 2\ \dist(\psi(0), \partial \psi(\mathbb{D}))\Phi_1(r).
$$
By assumption $|g'(z)|\leq k |h'(z)|$ in $\ID,$ where $k \in[0,1]$ and so, by Lemma~C 
and the classical Schwarz inequality, it follows that
\beqq		
\sum_{n=1}^\infty |b_n| \varphi_{n}(r) &\leq&  \sqrt{\sum_{n=1}^\infty |b_n|^2 \varphi_{n}(r)} \sqrt{\sum_{n=1}^{\infty}\varphi_{n}(r)}\\
&\leq& k \sqrt{\sum_{n=1}^\infty |a_n|^2 \varphi_{n}(r)} \sqrt{\sum_{n=1}^{\infty}\varphi_{n}(r)}\\
&\leq& 
2k\ \dist(\psi(0), \partial \psi(\mathbb{D}))\Phi_1(r).
\eeqq
Thus, we have
\beqq
\sum_{n=1}^\infty |a_n| \varphi_{n}(r)+\sum_{n=1}^\infty |b_n| \varphi_{n}(r) &\leq& 2(1+k)\dist(\psi(0), \partial \psi(\mathbb{D})) \Phi_1(r).\\
&\leq& \dist(\psi(0), \partial \psi(\mathbb{D})),\ \text{by Eqn.  \eqref{eqn 8}},
\eeqq
	for all $r \leq R,$ by the definition of $R.$ This proves the desired inequality \eqref{eqn 9}.
Now let us prove that $R$ is an optimal number. We consider the function
$$ \psi(z)=h(z) =\frac{1}{1-z}=\sum_{n=0}^{\infty}z^{n}, \ \ z \in \mathbb{D}
$$	
and $g'(z)= \lambda k h'(z),$ where $|\lambda|=1$. Then it is easy to see that
$$\dist(\psi(0), \partial \psi(\mathbb{D}))=\frac{1}{2}\ \,
\mbox{and}~g(z)=k \lambda \frac{z}{1-z}= k \lambda \sum_{n=1}^{\infty}z^{n},
$$
so that
\beqq
\sum_{n=1}^\infty |a_n| \varphi_{n}(r)+\sum_{n=1}^\infty |b_n| \varphi_{n}(r)=(1+k) \sum_{n=1}^{\infty}\varphi_{n}(r).
\eeqq
 Now it is easy to see that number is $> \frac{1}{2}$ when $r>R$. The proof of the theorem is complete.
\end{pf}

\beg
 Theorem \ref{thm 5} for the case of  $\varphi_{n}(r)=r^{n} \, (n \geq 1),$  gives the following result which was originally obtained at first in \cite[Theorem 1]{LiuPo}:
\beqq
\sum_{n=1}^\infty |a_n|r^n+\sum_{n=1}^\infty |b_n|r^n \leq \dist(\psi(0), \partial \psi(\mathbb{D})) \  \ \text{for}\  \ r \le \frac{1}{3+2k}.
\eeqq
The constant $\frac{1}{3+2k}$ is sharp.	
\eeg

\medskip
\noindent{\bf Acknowledgment.} I would like to thank my supervisor Prof. S. Ponnusamy for his support during the course of this work, fruitful discussions and valuable comments on this manuscipt.


\begin{thebibliography}{1}

\bibitem{Abu-1} Y. Abu-Muhanna,
Bohr phenomenon in subordination and bounded harmonic classes,
\emph{Complex Var. Elliptic Equ.} \textbf{55} (11)(2010), 1071--1078.
https://doi.org/10.1080/17476931003628190

\bibitem{AAPon1}  Y. Abu-Muhanna, R. M. Ali  and S. Ponnusamy,
On the Bohr inequality, In ``Progress in Approximation Theory and Applicable Complex Analysis'' (Edited by N.K. Govil et al. ),
Springer Optimization and Its Applications \textbf{117} (2016), 265--295.


\bibitem{AizenbergVidras} L. Aizenberg and A. Vidras,
On the Bohr radius for two classes of holomorphic functions,
\emph{Sibirsk. Mat. Zh.}, \textbf{45}(4) (2004), 734--746 (In Russian).
(English version: \emph{Sib. Math. J.}, \textbf{45}(4) (2004), 606--617).
https://doi.org/10.1023/B:SIMJ.0000035827.35563.b6.

\bibitem{Alkha}  S. A. Alkhaleefah, Bohr phenomenon for special family of analytic functions and harmonic mappings,
\emph{Probl. Anal. Issues Anal.} \textbf{9}(27)(3) (2020), 3--13. DOI: 10.15393/j3.art.2020.7990

\bibitem{AlkKayPon1} S.A. Alkhaleefah, I. R. Kayumov and S. Ponnusamy, On the Bohr inequality with a fixed zero coefficient,
\emph{Proc. Amer. Math. Soc.} \textbf{147}(12) (2019), 5263--5274.
https://doi.org/10.1090/proc/14634

\bibitem{BenDahKha} C. B\'{e}n\'{e}teau, A. Dahlner and D. Khavinson,
Remarks on the Bohr phenomenon,
\emph{Comput. Methods Funct. Theory} \textbf{4}(1) (2004), 1--19.
https://doi.org/10.1007/BF03321051

\bibitem{BoasKhavin-97-4} H. P. Boas and  D. Khavinson,
Bohr's power series theorem in several variables,
\emph{Proc. Amer. Math. Soc.} \textbf{125}(10) (1997),  2975--2979.
https://doi.org/10.1090/S0002-9939-97-04270-6

\bibitem{Bohr} H. Bohr, A theorem concerning power series,
\emph{Proc. London Math. Soc.} \textbf{13}(2) (1914), 1--5.
https://doi.org/10.1112/plms/s2-13.1.1

\bibitem{Bomb-1962} E. Bombieri,
Sopra un teorema di H. Bohr e G. Ricci sulle funzioni maggioranti delle serie di potenze,
\emph{Boll. Unione Mat. Ital.} \textbf{17} (1962), 276--282.

\bibitem{BombBour-2004}
E. Bombieri and J. Bourgain, A remark on Bohr's inequality,
\emph{IMRN International Mathematics Research Notices}, \textbf{80} (2004), 4307--4330.
https://doi.org/10.1155/S1073792804143444

\bibitem{DjaRaman-2000} P.~B.~Djakov and M.~S.~Ramanujan,
A remark on Bohr's theorems and its generalizations,
\emph{J. Analysis} \textbf{8} (2000), 65--77.

\bibitem{Duren} P. Duren, Univalent functions, Springer-Verlag, New York, 1983.

\bibitem{EvPoRa} S. Evdoridis, S. Ponnusamy and  A. Rasila,
Improved Bohr's inequality for locally univalent harmonic mappings,
\emph{Indag. Math. (N.S.)}, \textbf{30} (2019), 201--213.
https://doi.org/10.1016/j.indag.2018.09.008

\bibitem{GarMasRoss-2018}
S.~R.~Garcia, J.~ Mashreghi and W.~T.~Ross,
\textit{Finite Blaschke products and their connections}, Springer, Cham, 2018.

\bibitem{HHK2009}    H. Hamada, T. Honda, and G. Kohr,
Bohr's theorem for holomorphic mappings with values in homogeneous balls,
\emph{Israel J. Math.} \textbf{173} (2009), 177--187.
https://doi.org/10.1007/s11856-009-0087-9

\bibitem{KaykhaPo}  I. R. Kayumov, D. M. Khammatova and S. Ponnusamy,
Bohr inequality for the generalized  Ces\`{a}ro averaging operators, Preprint.

\bibitem{KayPon1} I. R. Kayumov and S. Ponnusamy,
Bohr inequality for odd analytic functions,
\emph{Comput. Methods Funct. Theory} \textbf{17} (2017), 679--688.
https://doi.org/10.1007/s40315-017-0206-2

\bibitem{KayPon3} I. R. Kayumov and S. Ponnusamy, Improved version of Bohr's inequality,
\emph{C. R. Math. Acad. Sci. Paris} \textbf{356}(3) (2018),  272--277. https://doi.org/10.1016/j.crma.2018.01.010

\bibitem{KayPon2} I. R. Kayumov and S. Ponnusamy,
Bohr's inequalities for the analytic functions with lacunary series and harmonic functions,
\emph{J. Math. Anal. and Appl.,}  \textbf{465} (2018), 857--871. https://doi.org/10.1016/j.jmaa.2018.05.038

\bibitem{KPS} I. R. Kayumov, S. Ponnusamy and N. Shakirov,
Bohr radius for locally univalent harmonic mappings,
\emph{Math. Nachr.}  \textbf{291} (2018), 1757--1768.
https://doi.org/10.1002/mana.201700068

\bibitem{LiuPo} Z. H. Liu, and S. Ponnusamy,
Bohr radius for subordination and K-quasiconformal harmonic mappings,
\emph{Bull. Malays. Math. Sci. Soc.} \textbf{42} (2019), 2151--2168.

https://doi.org/10.1007/s40840-019-00795-9

\bibitem{PoViWirth-3} S. Ponnusamy, R. Vijayakumar and K.-J. Wirths,
Modifications of Bohr's inequality in various settings, Preprint.
	
\end{thebibliography}
\end{document}